\documentclass[12pt]{amsart}
\usepackage{amsfonts,amssymb,amsthm,eucal,amsmath,verbatim,bbold}
\allowdisplaybreaks[3]
\newtheorem{thm}{Theorem}

\newcommand{\R}{\mathbb{R}}

\newcommand{\E}{\mathbb{E}}
\newcommand{\N}{\mathbb{N}}

\newcommand{\s}{\mathbb{S}}

\renewcommand{\O}{\mathcal{O}}

\newcommand{\vol}{\mathop{\mathrm{vol}}}

\newcommand{\inprod}[2]{\left\langle #1, #2 \right\rangle}

\newcommand{\ds}{\displaystyle}
\newcommand{\F}{\mathcal{F}}
\newcommand{\X}{\mathcal{X}}
\renewcommand{\P}{\mathbb{P}}

\renewcommand{\L}{\mathcal{L}}

\newcommand{\tr}{\mathrm{Tr\,}}

\renewcommand{\O}{\mathcal{O}}

\newcommand{\n}{\mathfrak{N}}

\begin{document}

\title{Approximation of projections of random vectors}
\author{Elizabeth Meckes}
\thanks{Research supported by an American Institute of Mathematics 
Five-year Fellowship and NSF grant DMS-0852898.}

\begin{abstract}
Let $X$ be a $d$-dimensional random vector and $X_\theta$ its projection
onto the span of a set of orthonormal vectors $\{\theta_1,\ldots,\theta_k\}$.
Conditions on the distribution of $X$ are given such that if $\theta$
is chosen according to Haar measure on the Stiefel manifold, the 
bounded-Lipschitz distance from $X_\theta$ to a Gaussian distribution 
is concentrated at its expectation; furthermore, an explicit bound is given for 
the expected distance, in terms of $d$, $k$, and the distribution of $X$,
allowing consideration not just of fixed $k$ but of $k$ growing with $d$.
The results are applied in the setting of projection pursuit, showing 
that most $k$-dimensional projections of $n$ data points in $\R^d$ are 
close to Gaussian, when $n$ and $d$ are large and $k=c\sqrt{\log(d)}$ for
a small constant $c$.

\end{abstract}

\maketitle

\section{Introduction}
There is a large class of results dealing with random variables (or
measures) defined in terms of a parameter (say, a point on the
sphere), which say that for a large measure of these parameters, the
behavior of the random variable is well-approximated by some model
distribution.  Early work in this direction was done by Sudakov \cite{Sud},
who showed that under some relatively mild conditions, most one-dimensional
marginals of a high-dimensional measure are close to each other.  This
line of research was further developed by von Weisz\"acker \cite{vonW},
who showed that the canonical distribution around which one-dimensional
marginals tend to cluster is close to a mixture of Gaussian distributions.
In both \cite{Sud} and \cite{vonW}, the results are about the limiting
behavior of one-dimensional projections, as the ambient dimension tends
to infinity, although von Weisz\"acker points out that one could extend
the methods to deal with higher fixed-dimensional projections, as the 
ambient dimension tends to infinity.  More recent work in this area was
done by Bobkov \cite{Bob}, who obtained concentration results
for the distance from a one-dimensional projection of an isotropic 
log-concave random vector to a Gaussian distribution.  

The purpose of this paper is to prove multivariate versions of
such theorems; that is, to consider rank $k$ projections of random vectors,
instead of just rank one.  Moreover, the approach yields results of 
a sufficiently quantitative nature to allow not only $k$ fixed, but $k$ 
growing with the ambient dimension.  The general case of approximating random
$k$-dimensional projections of probability measures on $\R^d$ is
considered, and is 
illustrated with an application to graphical projection
pursuit.  In particular, it is shown that typical $k$-dimensional
projections of $n$ data points in $\R^d$ are close to Gaussian for $n$
and $d$ large; the precise quantitative nature of the results yields 
limit theorems even for $k=c\log(d)$ for a small constant $c$. 
This
result generalizes the following univariate limit result of 
Diaconis and
Freedman.

\begin{thm}[Diaconis-Freedman \cite{dia-free}]\label{limit}
Let $x_1,\ldots,x_n$ be deterministic vectors in $\R^d$.  Suppose that
$n$, $d$ and the $x_i$ depend on a hidden index $\nu$, so that as $\nu$
tends to infinity, so do $n$ and $d$.  Suppose that there is a $\sigma^2>0$
such that, for all $\epsilon>0$,
\begin{equation}\label{limlengths}
\frac{1}{n}\Big|\left\{j\le n:\big||x_j|^2-\sigma^2d\big|>\epsilon d\right\}
\Big|\xrightarrow{\nu\to\infty}0,
\end{equation}
and suppose that
\begin{equation}\label{limorths}
\frac{1}{n^2}\Big|\left\{j,k\le n:\big|\inprod{x_j}{x_k}\big|>\epsilon d\right\}
\Big|\xrightarrow{\nu\to\infty}0.
\end{equation}
Let $\theta\in\s^{d-1}$ be distributed uniformly on the sphere, and consider 
the random measure $\mu_\nu^\theta$ which puts mass $\frac{1}{n}$ at 
each of the points $\inprod{\theta}{x_1},\ldots,\inprod{\theta}{x_n}$.  Then
as $\nu$ tends to infinity, the measures $\mu_\nu^\theta$ tend to $\n(0,\sigma^2)$
 weakly in probability.
\end{thm}

The method of proof here is described in a fairly specific context:
random measures indexed by points in the Stiefel manifold (one
could equivalently take points in the Grassman manifold), approximated
by Gaussian distributions.  However, the approach is quite general
and could in principle be adapted to a family of random measures indexed by
points in a metric
probability space possessing the concentration of measure phenomenon.  
Further, one could easily adapt the program
to deal with non-Gaussian limits.  In particular, Stein's method has been
used to prove approximation results for many other limiting distributions, 
e.g. Poisson \cite{chen75, agg89, agg90}; gamma 
\cite{luk94}; chi-square \cite{pick}; uniform  
on the discrete circle \cite{diaconis04}; the semi-circle law  
\cite{gottik05}; the binomial and multinomial distributions  
\cite{holmes04,loh92}; and the hypergeometric distribution \cite{holmes04};
these approaches could be combined with what is done here in order to
approximate by non-Gaussian distributions.

Before outlining the approach, some notation is needed.  The Euclidean
length of a vector $x\in\R^k$ is denoted $|x|$.  For an $n\times n$
matrix $M=\big[m_{ij}\big]_{i,j=1}^n$, the Hilbert-Schmidt norm is 
defined by 
$$\|M\|_{HS}=\tr(MM^T)=\sqrt{\sum_{i,j}m_{ij}^2}.$$
The Wasserstein distance between two random vectors $X$ and $Y$ is 
defined by 
$$d_W(X,Y)=\sup_{\{f:|f(x)-f(y)|\le|x-y|\}}\big|\E f(Y)-\E f(X)\big|.$$
The 
bounded-Lipschitz distance is defined by $$d_{BL}(X,Y):=\sup_{\|f\|_1\le 1}
\big|\E f(S)-\E f(Y)\big|,$$ where $$\|f\|_1:=\max\left\{\|f\|_\infty,
\sup_{x\neq y}\frac{|f(x)-f(y)|}{|x-y|}\right\}.$$

The class 
of $m$-times continuously differentiable functions on 
$\X\subseteq\R^d$ is denoted $C^m(\X)$, and has 
a norm defined by
$$\|f\|_m:=\sup_{0\le k\le m}\sup_{x\in\X}\|D^k f(x)\|_{op}.$$
Here, $D^kf(x)$ denotes the symmetric $k$-linear form given in components 
by 
$$D^kf(x)(y_1,\ldots,y_k):=\sum_{i_1,\ldots,i_k}\frac{\partial^k f}{\partial
x_{i_1}\cdots \partial x_{i_k}}(x)y_1^{i_1}\cdots y_k^{i_k},$$
where $y_j=(y_j^1,\ldots,y_j^d)$.  For an intrinsic definition of $D^kf(x)$,
see Federer \cite{fed}.  The ball of radius $R$ in $C^m(\X)$ with respect
to $\|\cdot\|_m$ is denoted $C^m_R(\X)$.

The Stiefel manifold $\mathfrak{W}_{d,k}$ is defined by 
$$ \mathfrak{W}_{d,k}=\{\theta=(\theta_1,\ldots,\theta_k):\theta_i\in\R^d, 
\inprod{\theta_i}{\theta_j}=\delta_{ij} \,\forall\, 1\le i,j\le k\},$$ 
with metric $\rho\big(\theta,\theta^\prime\big)=\left[\sum_{j=1}^k|\theta_j-
\theta_j^\prime|^2\right]^{1/2}$.  There is a unique rotation-invariant 
probability measure (Haar measure) on $\mathfrak{W}_{d,k}$; one way to 
construct 
it is by choosing $\theta_1$ uniformly from $\s^{d-1}$, then $\theta_2$ uniformly
from the orthogonal complement of $\theta_1$ in $\s^{d-1}$, and so on.

Now, suppose that a family of random vectors $X_\theta$ in $\R^k$ is indexed 
by $\theta\in\mathfrak{W}_{d,k}$.  The following is an outline of an approach 
to show that most $X_\theta$ are approximately Gaussian.

{\em 1. Prove an approximation result for the average distribution.}
If $X_\theta$ is defined fairly explicitly in terms of $\theta$, one 
can first try to use the following abstract normal approximation theorem to 
show that the average distribution of the $X_\theta$ (averaged over $\theta$
distributed according to Haar measure on $\mathfrak{W}_{d,k}$) 
is close to Gaussian.

\begin{thm}[\cite{CM}]\label{cont}
Let $X$ be a random vector in $\R^k$ and for each $\epsilon >0$ let $X_\epsilon$ be a random vector such that $\L(X)=
\L(X_\epsilon)$, with the property that $\lim_{\epsilon\to0}X_\epsilon=X$ 
almost surely.  Let $Z$ be a standard normal random vector in $\R^k$.
Suppose there is a function $\lambda(\epsilon)$ and a random matrix
$F$ such that the following conditions hold.
\begin{enumerate}
\item $$\frac{1}{\lambda(\epsilon)}\E\left[(X_\epsilon-X)_i\big|X\right]
\xrightarrow[\epsilon\to0]{L_1}- X. $$\label{first-diff}
\item $$\frac{1}{2\lambda(\epsilon)}\E\left[(X_\epsilon-X)(X_\epsilon-X)^T|
X\right]\xrightarrow[\epsilon\to0]{L_1}\sigma^2I_k+\E\left[F\big|X\right].$$
\label{second-diff}
\item For each $\rho>0$, $$\lim_{\epsilon\to0}\frac{1}{\lambda(\epsilon)}
\E\left[\big|X_\epsilon-X\big|^2
\mathbb{1}(|X_\epsilon-X|^2>\rho)\right]=0.$$\label{lindeberg}
\end{enumerate}
Then
\begin{equation}\label{contbd}
d_W(X,\sigma Z)\le\frac{1}{\sigma}\E\|F\|_{H.S.}
\end{equation}
\end{thm}
It should be pointed out that while this theorem is sufficiently general
for the applications carried out here, there is a more general 
version (see \cite {RR} or \cite{M}) 
allowing for approximations by Gaussian distributions
with non-trivial covariance matrices.  Furthermore, condition 
\ref{first-diff} need only hold approximately; see \cite{CM}.

In order to apply this theorem, an auxiliary random variable $X_{\theta,
\epsilon}$ must be constructed.  A natural construction which makes use of the 
symmetry of $\mathfrak{W}_{d,k}$ is to let $\theta_\epsilon$ be a 
``small random rotation'' of $\theta$ (this is made explicit in the 
applications to follow).  Then $(\theta,\theta_\epsilon)$ is an 
exchangeable pair of random points of $\mathfrak{W}_{d,k}$ by the rotation
invariance of the distribution of $\theta$, and so the random variables
$(X_\theta,X_{\theta_\epsilon})$ are also exchangeable and thus have the same 
distribution.  Furthermore, as $\epsilon\to0$, $\theta_\epsilon\to\theta$
almost surely, and so if $X_\theta$ is a continuous function of $\theta$, 
it will be true that $X_{\theta_\epsilon}\xrightarrow{\epsilon\to0}X_\theta$
almost surely. 

{\em 2. Use the concentration of measure on $\mathfrak{W}_{d,k}$ 
to show that for some distance $d(\cdot,\cdot)$, $d(X_\theta,\sigma Z)$ 
is close to its mean.}
 It is shown
in \cite{MS} that for $c_1=\sqrt{\frac{\pi}{2}}$ and $c_2=\frac{1}{8}$, 
for any $F:\mathfrak{W}_{d,k}\to\R$ with median $M_F$ and 
modulus of continuity $\omega_F(\eta)$, 
\begin{equation}
\P\big[|F(\theta_1,\ldots,\theta_k)-M_F|>\omega_F(\eta)\big]<c_1e^{-c_2\eta^2
d}.\label{stiefel}\end{equation}
Here, $\P$ is the rotation-invariant probability measure on 
$\theta_1,\ldots,\theta_k$ described above.  The median $M_F$ is a median with
respect to this measure.

Again, if the random variable $X_\theta$ is a sufficiently regular 
function of $\theta$, this theorem can be applied to the function $F(\theta)=
d_{BL}(X_\theta,\sigma Z),$ where $d_{BL}(X_\theta,\sigma Z)$ is the
conditional
bounded-Lipschitz distance from $X_\theta$ to $\sigma Z$, given $\theta$.  
Standard arguments
allow the median $M_F$ to be replaced by the mean $\E F(\theta)$, with
only minor loss.

{\em 3. Use entropy methods to bound $\E d_{BL}(X_\theta,\sigma Z)$.}
Consider the stochastic process $Y_f:=\big|\E_X f(X_\theta)-\E f(X_\theta)\big|$ 
indexed by the class of functions $\{f:\|f\|_1\le 1\}$ (or by some 
sub-class), where $\E_X$ denotes expectation with respect to $X$ only; that
is, conditional expectation with respect to the distribution of $X,\theta$,
conditioned on $\theta$.  Thus the bounded-Lipschitz distance from $X_\theta$
(given $\theta$) to its average distribution can be viewed as the 
supremum of a stochastic process. 
The same approach used to prove a concentration
result for $d_{BL}(X_\theta,\sigma Z)$ can be used to show that
$Y_f$  satisfies a
 sub-Gaussian increment condition of the type
$$\P\left[|Y_f-Y_g|>\epsilon\right]\le c_1e^{-\frac{c_2\epsilon^2}{\|f-g\|_1^2}},$$
for some constants $c_1$ and $c_2$.  For such a process, Dudley's entropy
bound can be used to estimate its supremum.  
Specifically, Dudley showed
the following.
\begin{thm}[Dudley, \cite{Dud}]
Let $\{X_t\}_{t\in T}$ be a stochastic process indexed by a metric space $T$
with distance $d$.  Suppose that there is a constant $c$ such 
that $X_t$ satisfies the increment condition
$$\forall u,\quad \P\left[|X_t-X_s|\ge u\right]\le c\exp\left(-\frac{u^2}{
2d(s,t)^2}\right).$$
Then there is a constant $C$ such that 
$$\E \sup_{t\in T}X_t\le C\int_0^\infty\sqrt{\log N(T,d,\epsilon)}d\epsilon,$$
where $N(T,d,\epsilon)$ is the $\epsilon$-covering number of $T$ with 
respect to the distance $d$.
\end{thm}
One can apply this theorem not to the index set $\{f:\|f\|_1\le 1\}$ (which 
has infinite $\epsilon$-covering number with respect to $\|\cdot\|_1$ for
$\epsilon<2$), but to a more restricted indexing set $\F$ of test functions.  
One may then be able to obtain a bound on $\E d_{BL}(X_\theta,\sigma Z)$
by approximation of functions $f$ with $\|f\|_1\le 1$ by functions from 
$\F$, together with the approximation for the average distribution proved
in Step 1.

\section{Random Projections}\label{projections-sec}
In this section, the method outlined in the introduction is applied 
in the case that $X$ is a random vector in $\R^d$, $\theta=(\theta_1,\ldots,
\theta_k)\in\mathfrak{W}_{d,k}$, and $X_\theta$
is the projection of $X$ onto the span of $\theta$; that is
$$X_\theta:=\big(\inprod{X}{\theta_1},\ldots,\inprod{X}{\theta_k}\big).$$
If $\theta$ is chosen randomly from $\mathfrak{W}_{d,k}$ (according to the 
rotation-invariant probability measure described in the introduction),
then the distributions of the $X_\theta$ are a family of random 
measures on $\R^k$ indexed by $\theta$.

To apply the method of the introduction, consider the random variable $X_\theta$
defined above, in the case that $\theta$ is chosen at random and 
independent of $X$.  The following
results describe the behavior of $X_\theta$, both on average and 
conditioned on $\theta$.

\begin{thm}\label{meanX}
Let $X$ be a random vector in $\R^n$, with $\E X=0$, $\E\left[|X|^2
\right]=\sigma^2d$, and $\E\big||X|^2\sigma^{-2}-d\big|:=A<\infty$.  
If $\theta$ is a random point of $\mathfrak{W}_{d,k}$ and $X_\theta$
is defined above,
$$d_W(X_\theta,\sigma Z)\le \frac{\sigma\sqrt{k}(A+1)+\sigma k}{d-1}.$$
\end{thm}

\begin{thm}\label{concX}
Suppose that $B$ is defined by
$B:=\sup_{\xi\in\s^{d-1}}\E\inprod{X}{\xi}^2.$  
For $\theta\in\mathfrak{W}_{d,k}$, let
$$d_{BL}(X_\theta,\sigma
  Z)=\sup_{\|f\|_1\le 1}\left|\E\left[ f(\inprod{X}{\theta_1},\ldots,\inprod{X}{
\theta_k})\big|\theta\right]-
  \E f(\sigma Z_1,\ldots,\sigma Z_k)\right|;$$ that is,
  $d_{BL}(X_\theta,\sigma Z)$ is the conditional bounded-Lipschitz distance
  from $X_\theta$ to $\sigma Z$, conditioned on $\theta$.  Then
for $\epsilon>2\pi\sqrt{\frac{B}{d}}$, and $\theta$ a random point of
$\mathfrak{W}_{d,k}$
$$\P\left[\big|d_{BL}(X_\theta,\sigma Z)-\E d_{BL}(X_\theta,\sigma Z)\big|
>\epsilon\right]\le \sqrt{\frac{\pi}{2}}e^{-\frac{d\epsilon^2}{32 B}}.$$

\end{thm}

\begin{thm}\label{dist-bdX}
There is a constant $C>1$ such that 
$$\E d_{BL}(X_\theta,\sigma Z)\le \frac{C^kB}{d^{\frac{2}{9k+4}}}+
\frac{\sigma\sqrt{k}(A+1)+\sigma k}{d-1}.$$
\end{thm}

Observe that together, Theorems \ref{concX} and \ref{dist-bdX} show that for
$\epsilon\ge\frac{2C^kB}{d^{\frac{2}{9k+4}}}+\frac{2\sigma\sqrt{k}(A+1)+2
\sigma k}{d-1}$, 
$$\P\left[d_{BL}(X_\theta,\sigma Z)>\epsilon\right]\le\sqrt{\frac{\pi}{2}}
e^{-\frac{d\epsilon^2}{2^7B}}.$$
Note that the bound on the right tends to zero as $d\to\infty$ for any 
$\epsilon$ in this range.

\begin{proof}[Proof of Theorem \ref{meanX}]
Observe first that $\E X_\theta=0$ by symmetry and 
\begin{equation*}\begin{split}
\E (X_\theta)_i(X_\theta)_j=\E\inprod{\theta_i}{X}\inprod{\theta_j}{X}&=
\sum_{r,s=1}^d\E\left[\theta_{ir}\theta_{js}\right]
\E\left[X_{r}X_{s}\right]=\frac{\delta_{ij}}{d}\E\left[|X|^2\right]
=\delta_{ij}\sigma^2,\end{split}\end{equation*}
where the second-last equality follows from $\E\left[\theta_{ir}\theta_{js}
\right]=\frac{1}{d}\delta_{ij}\delta_{rs}.$

To apply the Theorem \ref{cont} to $X_\theta$, one first has to construct 
$X_{\theta,\epsilon}$.  Let $$A_\epsilon:=
\begin{bmatrix}\sqrt{1-\epsilon^2}&\epsilon\\-\epsilon&
\sqrt{1-\epsilon^2}\end{bmatrix}\oplus I_{d-2}=I_d+\begin{bmatrix}
-\frac{\epsilon^2}{2}+\delta&\epsilon\\-\epsilon&
-\frac{\epsilon^2}{2}+\delta\end{bmatrix}\oplus 0_{d-2},$$ 
where $\delta=O(\epsilon^4)$.
Let $U\in\O_d$ be a random orthogonal
matrix, independent of $X$, and define 
$X_{\theta,\epsilon}:=
\left(\inprod{UA_\epsilon U^T\theta_1}{X},\ldots,\inprod{UA_\epsilon U^T
\theta_k}{X}\right)$; the pair $(X_\theta,X_{\theta,\epsilon})$ is exchangeable by
the rotation invariance of the distribution of $\theta$,
and so $\L(X_\theta)=\L(X_{\theta,\epsilon})$.  

Let $K$ be the $d\times 2$ matrix given by the first two columns of $U$ and
let $C=\begin{bmatrix}0&1\\-1&0\end{bmatrix}$; define the matrix $Q=\big[q_{ij}
\big]_{i,j=1}^d=KCK^T.$
Then, writing $X_\theta=(X^\theta_1,\ldots,X^\theta_k)$ and $X_{\theta,\epsilon}
=(X^\theta_{\epsilon,1},\ldots,X^\theta_{
\epsilon,k}),$
\begin{equation*}\begin{split}
\E\left[X^\theta_{\epsilon,j}-X^\theta_j\big|X,\theta\right]&=\E\left[
\inprod{(UA_\epsilon U^T-I)\theta_j}{X}\big|X,\theta
\right]\\&=\epsilon\E\left[\inprod{Q\theta_j}{X}\big|X,\theta\right]-\frac{
\epsilon^2}{2}\E\left[\inprod{KK^T\theta_j}{X}\big|X,\theta\right]+
O(\epsilon^4).
\end{split}\end{equation*}
Recall that $Q$ and $K$ are determined by $U$ alone, and that $U$ is 
independent of $X,\theta$.
It is easy to show that $\E\big[Q\big]=
0_d$ and $\E\big[KK^T\big]=\frac{2}{d}I_d$, thus
$$\E\left[X_{\theta,\epsilon}-X_\theta\big|X,\theta\right]=-\frac{\epsilon^2}{
d}X_\theta+O(\epsilon^4).$$
Condition \ref{first-diff} of Theorem \ref{cont} is thus satisfied with
$\F=\sigma(X,\theta)$ and $\lambda(\epsilon)=\frac{\epsilon^2}{d}$.

It is elementary but tedious to show that $\E q_{rs}q_{tv}=
\frac{2}{d(d-1)}\big[\delta_{rt}\delta_{sv}-\delta_{rv}\delta_{st}\big]$
(the computation is carried out in detail in \cite{CM}).  
Making use of this yields
\begin{equation*}\begin{split}
\E&\left[(X^\theta_{\epsilon,j}-X^\theta_j)(X^\theta_{\epsilon,\ell}-X^\theta_\ell
  )\big|X,\theta\right]\\& \qquad\qquad= \epsilon^2\E\left[\inprod{Q
    \theta_j}{X}\inprod{Q\theta_\ell}{X} \big|X,\theta\right]+
O(\epsilon^3)\\&\qquad\qquad=\epsilon^2\sum_{r,s,t,v=1}^d\E\left[q_{rs}q_{tv}
  \theta_{js}\theta_{\ell
    v}X_rX_{t}\big|X,\Theta\right]+O(\epsilon^3)\\&
\qquad\qquad=\frac{2\epsilon^2}{d(d-1)}\left[
  \sum_{r,s=1}^d\theta_{js}\theta_{\ell
    s}X_{r}^2-\sum_{r,s=1}^d\theta_{js} \theta_{\ell
    r}X_{r}X_{s}\right]+O(\epsilon^3)\\&\qquad\qquad
=\frac{2\epsilon^2}{d(d-1)}\left[\delta_{j\ell}|X|^2-X^\theta_j
  X^\theta_\ell\right]+
O(\epsilon^3)\\&\qquad\qquad=\frac{2\epsilon^2\sigma^2}{d}\delta_{j\ell}+\frac{2
  \epsilon^2}{d(d-1)}\left[\delta_{j\ell}\big(|X|^2-\sigma^2d\big)+\delta_{j\ell}
  \sigma^2-X^\theta_jX^\theta_\ell\right]+O(\epsilon^3).
\end{split}\end{equation*}
The random matrix $F$ of Theorem \ref{cont} 
is thus defined by 
$$F=\frac{1}{d-1}\left[\big(|X|^2-\sigma^2d\big)I_k+\sigma^2I_k-X_\theta 
X_\theta^T\right].$$
It follows from the theorem that 
\begin{equation}\begin{split}\label{annealX}
d_{W}(W,\sigma Z)&\le\frac{1}{\sigma}\E\|F\|_{H.S.}\\&\le\frac{\sigma
\sqrt{k}}{d-1}\left[\E\left|\frac{
|X|^2}{\sigma^2}-d\right|+1\right]+\frac{\sigma}{d-1}\E\left[\sum_{j}
\left(\frac{X^\theta_j}{\sigma}\right)^2\right]\\&\le\frac{\sigma\sqrt{k}(A+1)+
\sigma k}{d-1}.\end{split}\end{equation}

\end{proof}

\medskip

\begin{proof}[Proof of Theorem \ref{concX}]
Define a function $F:\mathfrak{W}_{d,k}\to\R$ by 
$$F(\theta)=\sup_{\|f\|_1\le 1}\big|\E_X f(X_\theta))-\E f(\sigma
Z)\big|,$$ where $\E_X$ denotes the expectation with respect to the 
distribution of $X$ only;
that is, 
$$\E_X f(X_\theta)=\E\left[f(X_\theta)\big|\theta\right].$$

To apply the concentration of measure on $\mathfrak{W}_{d,k}$, it
is necessary to determine the modulus of continuity of $F$.  
First, observe that for $f$ with $\|f\|_1\le 1$ given,
\begin{eqnarray*}
\Big|\big|\E_X f(X_\theta)-\E f(\sigma Z)\big|&-&\big|\E_X
f(X_\theta')- \E f(\sigma Z)\big|\Big|\\&\le&\Big|\E_X
f(X_\theta')-\E_X f(X_\theta))\Big|\\&=&
\E\left[f\big(\inprod{X}{\theta_1'},\ldots,\inprod{X}{\theta_k'}\big)-
f\big(\inprod{X}{\theta_1},\ldots,\inprod{X}{\theta_k}\big)\Big|\theta,\theta'
\right]\\&\le&\E\left[\big|\big(\inprod{X}{\theta_1'-\theta_1},\ldots,
\inprod{X}{\theta_k'-\theta_k}\big)\big|\Big|\theta,\theta'\right]
\\&\le&\sqrt{\sum_{j=1}^k|\theta_j'-\theta_j|^2\E\inprod{X}{\frac{
\theta_j'-\theta_j}{|\theta_j'-\theta_j|}}^2}\\&\le&
\rho(\theta,\theta')\sqrt{B}.
\end{eqnarray*}
It follows that
\begin{equation*}\begin{split}
\Big|d_{BL}(X_\theta,\sigma Z)&-d_{BL}(X_{\theta'},\sigma Z)\Big|\\&=
\left|\sup_{\|f\|_1\le 1}\big|\E_X f(X_\theta)-\E f(\sigma Z)\big|-
\sup_{\|f\|_1\le 1}\big|\E_X f(X_{\theta'})-\E f(\sigma Z)\big|\right|\\&\le
\sup_{\|f\|_1\le 1}\Big|\big|\E_X f(X_\theta)-\E f(\sigma Z)\big|-\big|
\E_X f(X_{\theta'})-\E f(\sigma Z)\big|\Big|\\&\le
\rho(\theta,\theta')\sqrt{B},
\end{split}\end{equation*}
thus $d_{BL}(X_\theta,\sigma Z)$ is a Lipschitz function on 
$\mathfrak{W}_{k,d}$, with Lipschitz constant $\sqrt{B}$.
Applying the concentration of measure inequality from Inequality \eqref{stiefel}
of the introduction then implies
that 
$$\P\big[|F(\theta_1,\ldots,\theta_k)-M_F|>\epsilon\big]<\sqrt{\frac{
\pi}{2}}e^{-\frac{\epsilon^2d}{8B}}.$$

Now, if $\theta=(\theta_1,\ldots,\theta_k)$ is a Haar-distributed 
random point of $\mathfrak{W}_{d,k}$, then
\begin{equation*}\begin{split}
\big|\E F(\theta)-M_F\big|\le\E\big|F(\theta)-M_F\big|=
\int_0^\infty\P\Big[\big|
F&(\theta)-M_F\big|>t\Big]dt\\&\le\int_0^\infty \sqrt{\frac{\pi}{2}}e^{-
\frac{dt^2}{8B}}dt
=\pi\sqrt{\frac{B}{d}}.
\end{split}\end{equation*}
So as long as $\epsilon>2\pi\sqrt{\frac{B}{d}},$ 
replacing the median of $F$ with its mean only changes the constants:
\begin{equation*}\begin{split}
\P\left[\left|F(\theta)-\E F(\theta)\right|>\epsilon\right]\le
\P\Big[\big|F(\theta)-M_F&\big|>\epsilon-\big|M_F-\E F(\theta)
\big|\Big]\\&\le\P\left[\left|F(\theta)
-M_F\right|>\frac{\epsilon}{2}\right]\le \sqrt{\frac{\pi}{2}}
e^{-\frac{d\epsilon^2}{32B}}.
\end{split}\end{equation*}
\end{proof}

What has just been shown is that $d_{BL}(X_\theta,\sigma Z)$ is concentrated 
about its mean; it remains to give a bound for this mean (Theorem 
\ref{dist-bdX}).

\begin{proof}[Proof of Theorem \ref{dist-bdX}]
As indicated in the introduction, Theorem \ref{dist-bdX} is proved 
making use of Dudley's entropy bound for bounding the expected
value of the supremum of a stochastic process.  Let $X_f:=\big|\E_Xf(X_\theta)-
\E f(X_\theta)\big|$.  Then $\{X_f\}_f$ is a stochastic process (each $X_f$
is a random variable depending on $\theta$) indexed by a family of
functions $f$.  The same type of concentration argument used above can
be used to show that this process is sub-Gaussian.  

Let $f:\R^k\to\R$ be Lipschitz with Lipschitz constant $L$ and
 consider the 
 function $G=G_f$ defined on $\mathfrak{W}_{d,k}$ by 
 $$G(\theta_1,\ldots,\theta_k)=\E_X f(X_\theta)=
\E\left[f(\inprod{\theta_1}{X},\ldots,\inprod{\theta_k}{X})\big|\theta\right].$$
The same argument as above shows that $G$ is Lipschitz on 
$\s^{d-1}$ with Lipschitz constant $L\sqrt{B}$.
 It thus follows from \eqref{stiefel} that 
$$\P\left[\left|G(\theta)-M_G\right|>\epsilon\right]\le \sqrt{\frac{\pi}{2}}
e^{-\frac{d\eta^2}{8L^2B}},$$and if $\epsilon>2\pi\sqrt{\frac{L^2B}{d}},$ 
\begin{equation}\begin{split}\label{GconcX}
\P\left[\left|G(\theta)-\E G(\theta)\right|>\epsilon\right]
&\le\sqrt{\frac{\pi}{2}}e^{-\frac{d\epsilon^2}{32L^2B}}.
\end{split}\end{equation}
Observe that, for $\theta$ a Haar-distributed random point of $\mathfrak{W}_{
d,k}$, $\E G(\theta)=\E f(X_\theta)$,
and so \eqref{GconcX} can be restated as
$$\P\left[X_f>\epsilon\right]\le
\sqrt{\frac{\pi}{2}}\exp\left[-\frac{d\epsilon^2}{2^7L^2B}\right].$$

Note that 
\begin{equation*}\begin{split}
\big|X_f-X_g\big|=\Big|\big|\E_Xf(X_\theta)-\E& f(X_\theta)\big|-
\big|\E_Xg(X_\theta)-\E g(X_\theta)\big|\Big|\\&
\le\Big|\E_X(f-g)(X_\theta)-\E
(f-g)(X_\theta)\Big|=X_{f-g},\end{split}\end{equation*}
thus for $\epsilon>4\pi L(f-g)\sqrt{\frac{B}{d}}$, for $L(f-g)$ the Lipschitz
constant of $f-g$,
$$\P\left[\big|X_f-X_g\big|>\epsilon\right]\le\P\left[X_{f-g}>\epsilon
\right]\le\sqrt{\frac{\pi}{2}}\exp\left[\frac{-d\epsilon^2}{2^7[L(f-g)]^2B}
\right]\le\sqrt{\frac{\pi}{2}}\exp\left[\frac{-d\epsilon^2}{2^7\|f-g\|_1^2B}
\right].$$
The condition on $\epsilon$ may be removed by replacing the factor
of $\sqrt{\frac{\pi}{2}}$ in the bound above by, e.g., $3\sqrt{\frac{\pi}{2}}$.
The process $\{X_f\}$
therefore satisfies the sub-Gaussian increment condition for the distance 
$d^*(f,g):=\frac{8\sqrt{B}}{\sqrt{d}}\|f-g\|_1.$

Consider the class $C^m_1(B_R)$ of functions $f$ which are supported on $B_R:=
\{x\in\R^k:|x|\le R\}$ such that $\|f\|_m:=\sup_{0\le j\le m}\sup_{x\in B_R}\|
D^jf(x)\|_{op}\le 1.$  It is proved in the appendix that for 
$\epsilon<2$ and $m\ge 2$, the $\epsilon$-covering
number for this set with respect to the norm $\|\cdot\|_1$ is bounded
by 
$$\exp\left[\left((3\log(5)-\frac{m}{m-1}\log(\epsilon)+
\frac{c_1}{\epsilon^{\frac{k}{m-1}}}\right)e^{k+m-2}\right]$$
with
$$c_1=\frac{2\pi^{k/2}(R+1)^k\big((m+4)\log(2)\big)\big(5\big)^{\frac{k}{m-1}}}{
k\Gamma\left(\frac{k}{2}\right)}.$$
It follows that the $\epsilon$-covering number with respect to the 
distance $d^*$ is bounded by
\begin{equation*}\begin{split}
\exp&\left[\left(3\log(5)-\frac{m}{(m-1)}\log\left(\frac{\epsilon\sqrt{d}}{
8\sqrt{B}}\right)\phantom{\frac{\left[40\sqrt{B}\right]^{\frac{k}{m
-1}}}{\left[\epsilon\sqrt{d}\right]^{\frac{k}{m-1}}}}\right.\right.
\\&\qquad\qquad\left.\left.+\frac{2\log(2)(m+4)\left[\sqrt{\pi}(R+1)\right]^k
\left[40\sqrt{B}\right]^{\frac{k}{m-1}}}{k\Gamma\left(\frac{k}{2}\right)\left[
\epsilon\sqrt{d}\right]^{\frac{k}{m-1}}}\right)e^{k+m-2}\right].
\end{split}\end{equation*}

Since functions $f\in C^m_1(B_R)$ have in particular $\|f\|_1\le1$, 
this class also satisfies the sub-Gaussian increment condition with
respect to the metric $d^*$.  Note that the diameter of $C^m_1(B_R)$
with respect to $d^*$ is bounded above by $16\sqrt{\frac{B}{d}}.$  It follows
from Dudley's entropy bound that there is a constant $C$ such that 
$\ds\E \Big[\sup_{f\in C^m_1(B_R)}X_f\Big]$ is bounded above by 
$$Ce^{\frac{k+m-2}{2}}\int_0^{16\sqrt{\frac{B}{d}}}\sqrt{3\log(5)-\frac{m}{(m-1)}
\log\left(\frac{\epsilon\sqrt{d}}{8\sqrt{B}}\right)+\frac{2\log(2)(m+4)
\left[\sqrt{\pi}(R+1)\right]^k\left[40\sqrt{B}\right]^{\frac{k}{m-1}}}{
k\Gamma\left(\frac{k}{2}\right)\left[\epsilon\sqrt{d}\right]^{\frac{k}{m-1}}}}
d\epsilon.$$
Making the substitution $s=\frac{\epsilon\sqrt{d}}{8\sqrt{B}}$ then gives
an upper bound of 
$$C'e^{\frac{k+m-2}{2}}\sqrt{\frac{B}{d}}\int_0^2\sqrt{3\log(5)-\frac{m}{m-1}
\log(s)+\frac{2\log(2)(m+4)
\left[\sqrt{\pi}(R+1)\right]^k(5)^{\frac{k}{m-1}}}{
k\Gamma\left(\frac{k}{2}\right)s^{\frac{k}{m-1}}}}ds$$ for another 
constant $C'$.
Looking at the first two summands and the third separately, 
as long as $m>\frac{k}{2}+1$, 
this implies that there is an absolute constant
$C$ such that $\ds\E \Big[\sup_{f\in C^m_1(B_R)}X_f\Big]$ is bounded by
$$\sqrt{\frac{B}{d}}\left(\frac{C^{k+m}R^{k/2}m^{3/2}}{(2m-k-2)
\sqrt{k\Gamma\left(\frac{k}{2}\right)}}\right),$$
or, as will be needed in what follows,
\begin{equation}\label{bdsmoothX}
\E\left[\sup_{f\in C^m_M(B_R)}X_f\right]\le
\sqrt{\frac{B}{d}}\left(\frac{MC^{k+m}R^{k/2}
  m^{3/2}}{(2m-k-2)\sqrt{k\Gamma\left(\frac{k}{2}\right)}}\right).
\end{equation}
From this bound, one can obtain a bound on $\E d_{BL}(X_\theta,)$
as follows.  Let 
$$\varphi_R(x)=\begin{cases}1&|x|\le R,\\R+1-|x|&R\le|x|\le R+1,\\0&R+1\le|x|;
\end{cases}$$
that is, $\varphi_R$ is a radially symmetric cut-off function with 
$\|\varphi_R\|_1\le1$, supported on $B_{R+1}$ and with $\varphi_R\equiv
1$ on $B_R$. For $f\in C^1_1(\R^k)$, let $f_R:=f\cdot\varphi_R$.  
Then $$\|f_R\|_1=\max\left\{\sup_x|f(x)\varphi_R(x)|,\sup_x\left|f(x)
\cdot\nabla\varphi_R(x)+\varphi_R(x)\nabla f(x)\right|\right\}\le 2.$$
Since $|f(x)-f_R(x)|=0$ if $x\in B_R$ and $|f(x)-f_R(x)|\le 1$ for 
all $x\in\R^k$, 
$$\big|\E_X f(X_\theta)-\E_X f_R(X_\theta)\big|\le\P\big[|X_\theta|>R\big|
\theta\big]\le
\frac{1}{R^2}\sum_{i=1}^k\E\big[
\inprod{X}{\theta_i}^2\big]\le\frac{Bk}{R^2},$$
and the same holds if $\E_X$ is replaced by $\E$.
It follows that 
\begin{equation}\label{truncX}
\left|X_f-X_{f_R}\right|\le\frac{2Bk}{R^2}.\end{equation}

Next, let $\psi:\R\to\R$ be a $C^\infty$ bump function, such that $0\le \psi(y)
\le 1$ for all $y$, $\psi(y)=1$ for $-1\le y\le 1$, $\psi(y)=0$ for $|y|>2$, 
and such that 
\begin{equation}\label{1dderivbdsX}
\left|\frac{d^j\psi}{dy^j}(y)\right|\le C^j j^{2j}
\end{equation} for 
all $j\in\N$ (the existence of such a function is guarranteed by Theorem 1.4.2
of \cite{Hor}).   
For $x\in\R^k$, define 
$$\psi_t(x)=\frac{C(k)}{t^k}\psi\left(\frac{|x|}{t}\right),$$
where $C(k)$ is a constant depending only on $k$, such that 
$\int_{\R^k}\psi_t=1.$  
Observe that it follows from the bounds \eqref{1dderivbdsX} that
\begin{equation}\label{derivbdsX}
\|D^j\psi_t(x)\|_{op}\le\frac{C(k)C^jj^{2j}}{t^{k+j}}\mathbb{1}(t\le|x|\le 2t).
\end{equation}
For $g\in C^1_2(\R^k)$, let $g_t(x):=g\ast\psi_t(x).$  
Let $Y_t$ be a random vector in $\R^k$ with density $\psi_t$,
independent of $X,\theta$.  Then one 
can write
$$\E_X g_t(X_\theta)=\E_X g(X_\theta+Y_t),$$
and the same with $\E$ in place of $\E_X$.
Since $g\in C^1_2(\R^k)$, it follows that 
$$\max\Big(\big|\E_X g(X_\theta)-\E_X g_t(X_\theta)\big|,\big|\E g(X_\theta)
-\E g_t(X_\theta)\big|\Big)\le 2\E|Y_t|\le 4t,$$
from which it follows that 
\begin{equation}\label{smoothingX}
\left|X_g-X_{g_t}\right|\le8t.\end{equation}
Furthermore,
by Young's inequality, for $j\le m$, 
\begin{equation*}\begin{split}
\|D^jg_t(x)\|_{op}&\le\|g\|_\infty\int_{\R^k}\|D^j\psi_t(y)\|_{op}dy\le
\frac{2C(k)C^jj^{2j}}{t^{k+j}}\vol(B_{2t})=\frac{2^{k+2}\pi^{k/2}C(k)C^jj^{2j}}{
t^{j}k\Gamma\left(\frac{k}{2}\right)}.\end{split}\end{equation*}
Now, integrating in polar coordinates,
\begin{equation*}\begin{split}
\frac{1}{C(k)}&=\int_{\R^k}t^{-k}\psi\left(\frac{|x|}{t}\right)dx
=\frac{2\pi^{k/2}}{\Gamma\left(\frac{k}{2}
\right)}\int_0^2\psi(r)r^{k-1}dr\ge\frac{2\pi^{k/2}}{\Gamma\left(\frac{k}{2}
\right)}\int_0^1r^{k-1}dr=\frac{2\pi^{k/2}}{k\Gamma\left(\frac{k}{2}
\right)}.\end{split}\end{equation*}
It follows that
$$\|D^jg_t(x)\|_{op}\le \frac{2^{k+1}C^jj^{2j}}{t^{j}}$$
for all $x\in\R^k$, and so $\|g_t\|_m\le\frac{2^{k+1}C^mm^{2m}}{t^{m}}$.
Finally, if $g$ is supported on $B_{R+1}$, then it is easy to 
see that $g_t$ is supported on $B_{R+1+2t}$.  

It now follows from \eqref{bdsmoothX}, \eqref{truncX} and \eqref{smoothingX}
that
\begin{equation}\begin{split}\label{finalX}
\E\left[\sup_{f\in C^1_1(\R^k)}X_f\right]&\le \E\left(\sup_{f\in
  C^1_1(\R^k)}\left[\left|X_f-X_{f_R}
  \right|+\left|X_{f_R}-X_{(f_R)_t}\right|+X_{(f_{R})_t}\right]\right)\\&\le
\frac{2Bk}{R^2}+8t+\sqrt{\frac{B}{d}}\left(\frac{2^{k+1}m^{2m}
  C^{k+m}(R+1+2t)^{k/2}m^{3/2}}{(2m-k-2)t^m\sqrt{k\Gamma\left(\frac{k}{2}
    \right)}}\right)\\&\le\frac{2Bk}{R^2}+8t+\sqrt{\frac{B}{d}}
\left(\frac{C^{k+m}m^{2m}R^{k/2}m^{3/2}}{(2m-k-2)t^m\sqrt{k\Gamma\left(\frac{k}{2}
    \right)}}\right).
\end{split}\end{equation}
Choosing $t=\frac{k}{R^2}$ yields
\begin{equation}\label{simpler1X}
\E\left[\sup_{f\in C^1_1(\R^k)}X_f\right]\le\frac{(2B+8)k}{R^2}+
\sqrt{\frac{B}{d}}\left(\frac{C^{k+m}m^{2m}R^{2m+k/2}m^{3/2}}{
(2m-k-2)k^m\sqrt{k\Gamma\left(\frac{k}{2}\right)}}\right).
\end{equation}
Now choosing $m=k$ and applying Stirling's formula to $\Gamma\left(
\frac{k}{2}\right)$ yields
\begin{equation}\label{simpler2X}
\E\left[\sup_{f\in C^1_1(\R^k)}X_f\right] \le\frac{(2B+8)k}{R^2}+
\sqrt{\frac{B}{d}}\left[\big(Ck^{3/4}\big)^{k}R^{9k/2}\right].
\end{equation}
Setting $R=\left(\frac{d}{k^{\frac{3k}{2}-2}}\right)^{\frac{1}{9k+4}}$
yields
\begin{equation}\label{simplestX}
\E\left[\sup_{f\in C^1_1(\R^k)}X_f\right]\le\frac{C^kB}{d^{\frac{2}{9k+4}}}.
\end{equation}
Finally, by Theorem \ref{meanX} and \eqref{simplestX}, 
\begin{equation*}\begin{split}
\E d_{BL}(X_\theta,\sigma Z)&\le\E\left(\sup_{f\in C^1_1(\R^k)}\Big[
|\E_Xf(X_\theta)-\E f(X_\theta)|+|\E f(X_\theta)-\E f(\sigma Z) |\Big]
\right)\\&\le \frac{C^kB}{d^{\frac{2}{9k+4}}}+\frac{\sigma\sqrt{k}(A+1)+\sigma k}{
d-1}.\end{split}\end{equation*}

\end{proof}

\section{Application: Projection Pursuit}\label{projpursuit-sec}
In this section, the theorems of the previous section are applied to
prove a quantitative, higher-dimensional version of a result of 
Diaconis and Freedman \cite{dia-free}.
Let $x_1,\ldots,x_n$ be deterministic vectors in $\R^d$; write 
$x_i=(x_{i,1},\ldots,x_{i,d}).$  Define $\sigma>0$ by the condition
$\frac{1}{n}\sum_{i=1}^n|x_i|^2=\sigma^2d$, and define $A$ and 
$B$ by $A:=\frac{1}{n}\sum_{i=1}^n\big|\sigma^{-2}|x_i|^2-d\big|$ and
$B:=\sup_{\theta\in\s^{d-1}}\frac{1}{n}\sum_{i=1}^n\inprod{\theta}{x_i}^2$.  
Observe that $\frac{\sigma^2d}{n}\le B\le \sigma^2 d$.  Also,if 
$X$ is distributed uniformly over the points $\{x_i\}$, then these
definitions of $\sigma$, $A$, and $B$ correspond to those in the 
previous section.

Let $\theta=(\theta_1,\ldots,\theta_k)$ be a random point in 
$\mathfrak{W}_{d,k}$, distributed according to the rotation-invariant 
probability measure described in the introduction, and 
consider the family of random measures $\mu_{n,d,k}^\theta$ defined in terms
of $\theta$ by 
$$\mu_{n,d,k}^\theta:=\frac{1}{n}\sum_{i=1}^n\delta_{(\inprod{\theta_1}{x_i},\ldots,
\inprod{\theta_k}{x_i})}.$$
That is, $\mu_{n,d,k}^\theta$ puts equal mass at the projections of each of 
the $x_i$ onto the span of $\theta_1,\ldots,\theta_k$. 

In Diaconis and Freedman \cite{dia-free} it was shown that, 
in the case $k=1$, the measures
$\mu_{n,d,1}^\theta$ converge weakly in probability to Gaussian as 
$n$ and $d$ tend to infinity, under the conditions that, for some $\sigma^2>0$
such that, for all $\epsilon>0$,
\begin{equation}\label{limlengths1}
\frac{1}{n}\Big|\left\{j\le n:\big||x_j|^2-\sigma^2d\big|>\epsilon d\right\}
\Big|\xrightarrow{\nu\to\infty}0,
\end{equation}
and 
\begin{equation}\label{limorths1}
\frac{1}{n^2}\Big|\left\{j,k\le n:\big|\inprod{x_j}{x_k}\big|>\epsilon d\right\}
\Big|\xrightarrow{\nu\to\infty}0.
\end{equation}
Here, $n$, $d$, and the $x_i$ depend on a hidden index $\nu$ such that as $\nu$
tends to infinity, so do $n$ and $d$.  A reasonable quantitative analog
would be to require $A$ and $B$ above to be bounded, independent of $n$ 
and $d$.  One could also allow them to grow slowly, as is clear from the 
statements of the theorems below.  Recall that $B\ge\frac{\sigma^2d}{n},$
so if $B$ is to remain bounded as $d$ tends to infinity, $n$ must tend
to infinity at least as fast as $d$.

In recent work of the author \cite{M2},
a quantitative version of the Diaconis-Freedman
 result was proved, giving an explicit bound
on $\P\left[d_{BL}(\mu_{n,d,1}^\theta,\gamma_{\sigma^2})\ge\epsilon\right],$
where $\gamma_{\sigma^2}$ is the Gaussian distribution on $\R$ with mean zero
and variance $\sigma^2$.  The results of Section \ref{projections-sec}
apply immediately to the random vector $X$ uniformly distributed 
on the $n$ points $\{x_i\}_{i=1}^n$ to give the following $k$-dimensional
extensions.

\begin{thm}\label{mean}
If $\theta$ is a random point of $\mathfrak{W}_{d,k}$ and $X_\theta$
is distributed according $\mu_{n,d,k}^\theta$, then
$$d_W(X_\theta,\sigma Z)\le \frac{\sigma\sqrt{k}(A+1)+\sigma k}{d-1}.$$
\end{thm}

\begin{thm}\label{conc}For $\theta\in\mathfrak{W}_{d,k}$, let
$$d_{BL}(X_\theta,\sigma
  Z)=\sup_{\|f\|_1\le 1}\left|\frac{1}{n}\sum_{i=1}^n
  f\left(\inprod{x_i}{\theta_1},\ldots,\inprod{x_i}{\theta_k}\right)-
  \E f(\sigma Z_1,\ldots,\sigma Z_k)\right|;$$ that is,
  $d_{BL}(X_\theta,\sigma Z)$ is the conditional bounded-Lipschitz distance
  from $X_\theta$ to $\sigma Z$, conditioned on $\theta$.  Then
for $\epsilon>2\pi\sqrt{\frac{B}{d}}$, and $\theta$ a random point of
$\mathfrak{W}_{d,k}$
$$\P\left[\big|d_{BL}(X_\theta,\sigma Z)-\E d_{BL}(W(\theta),\sigma Z)\big|
>\epsilon\right]\le \sqrt{\frac{\pi}{2}}e^{-\frac{d\epsilon^2}{32 B}}.$$

\end{thm}

\begin{thm}\label{dist-bd}
There is a constant $C>1$ such that 
$$\E d_{BL}(X_\theta,\sigma Z)\le \frac{C^kB}{d^{\frac{2}{9k+4}}}+
\frac{\sigma\sqrt{k}(A+1)+\sigma k}{d-1}.$$
\end{thm}

Observe that together, Theorems \ref{conc} and \ref{dist-bd} show that for
$\epsilon\ge\frac{2C^kB}{d^{\frac{2}{5k+4}}}+\frac{\sigma\sqrt{k}(A+1)+\sigma k}{d-1}$, 
$$\P\left[d_{BL}(X_\theta,\sigma Z)>\epsilon\right]\le\sqrt{\frac{\pi}{2}}
e^{-\frac{d\epsilon^2}{2^7B}}.$$
Note that the bound on the right tends to zero as $d\to\infty$ for any 
$\epsilon$ in this range.  In particular, if $A$ and $B$ are 
bounded and $\epsilon>0$
fixed, if $k=c\log(d)$, where $c$ is a sufficiently small constant (depending
on $\epsilon$), then 
$ \P\left[d_{BL}(X_\theta,\sigma Z)>\epsilon\right]$ decays
exponentially as $d$ tends to infinity.

\section{Appendix: The covering number of the class $C^m_M(\X)$}

Consider the class $C^m(\X)$
of $m$-times continuously differentiable functions on 
$\X\subseteq\R^d$ 
with norm defined by
$$\|f\|_m=\sup_{0\le k\le m}\sup_{x\in\X}\|D^k f(x)\|_{op}.$$
Here, $D^kf(x)$ denotes the symmetric $k$-linear form given in components 
by 
$$D^kf(x)(y_1,\ldots,y_k):=\sum_{i_1,\ldots,i_k}\frac{\partial^k f}{\partial
x_{i_1}\cdots \partial x_{i_k}}(x)y_1^{i_1}\cdots y_k^{i_k},$$
where $y_j=(y_j^1,\ldots,y_j^d)$.  For an intrinsic definition of $D^kf(x)$,
see Federer \cite{fed}.

Let $C^m_M(\X)$ be the ball of radius $M$ of $C^m(\X)$ with respect
to $\|\cdot\|_m$; in this section, the 
$\epsilon$-covering number of $C_1^m(\X)$ 
with respect to the norms $\|\cdot\|_{\infty}$
(defined the usual way; in our notation, this is $\|\cdot\|_0$)
and $\|\cdot\|_{1}$ is calculated for $m\ge 2$.  
  The proof closely follows 
the approach in van der Vaart and Wellner \cite{vanWel} but
uses the definition of $D^kf$ as a $k$-linear form instead of working in
coordinates with the partial derivatives of $f$.

First, choose a $\delta$-net $\{y_i\}_{i=1}^n$ of $\X$, with $\delta=\delta(
\epsilon)$ to be determined.  One can choose such a net so that $n\le\frac{
\vol(\X_1)}{\delta^d},$ where $\X_1:=\{x\in\R^d:\inf_{y\in\X}|x-y|\le1\}$.
Now, associate to each $f\in C^m_1(\X)$ an $(m-1)\times n$ array of
operators in the following way.  In the space of symmetric $k$-linear forms 
on $\R^d$, choose a $\frac{\delta^{m-k}}{2}$-net $\{T_i\}_{i=1}^M$, 
with respect to the operator norm.  The $(i,j)$-th entry of the array $A_f$
associated to $f$ is chosen to be the closest point in the appropriate
net to the $i$-linear form $D^if(y_j)$.   One can choose $\delta_0$ and 
$\delta_1$ such 
that if $f,g\in C^m_1(\X)$ have $A_f=A_g$ (with respect to either the $\delta_0$
 or the $\delta_1$ nets), then $\|f-g\|_{\infty}\le\epsilon$ for $\delta_0$
and $\|f-g\|_1\le\epsilon$ for $\delta_1$,
as follows.
For $x\in\X$ given, choose $y_i$ with $|x-y_i|\le\delta$.  By Taylor's theorem
applied to $f-g$, 
$$(f-g)(x)=\sum_{k=0}^{m-1}\frac{1}{k!}\inprod{D^k(f-g)(y_i)}{
(x-y_i,\ldots,x-y_i)}+R,$$
with $|R|\le\frac{2|x-y_i|^{m}}{m!}\le\frac{2\delta^{m}}{m!}$.
Since $A_f=A_g$, it follows that $\|D^k(f-g)(y_i)\|_{op}\le\delta^{m-k}$ 
for $1\le k\le m-1$, thus by the expansion above,
\begin{equation*}\begin{split}
|(f-g)(x)|&\le\sum_{k=0}^{m-1}\frac{\delta^k}{k!}\|D^k(f-g)(y_i)\|_{op}+
\frac{2\delta^{m}}{m!}\\&\le\sum_{k=0}^{m-1}\frac{\delta^k}{k!}\delta^{m-k}+
\frac{2\delta^{m}}{m!}\\&<5\delta^{m},
\end{split}\end{equation*}
since $m\ge 1$.
It follows that choosing  $\delta_0=\left(\frac{\epsilon}{5}\right)^{
\frac{1}{m}}$ means that if $A_f=A_g$ then $\|f-g\|_\infty<\epsilon.$

To choose $\delta_1$, apply Taylor's theorem to $D(f-g)$: if $|v|=1$,
$$\inprod{D(f-g)(x)}{v}=\sum_{k=1}^{m-1}\frac{1}{k!}\inprod{D^k(f-g)(y_i)}{
(x-y_i,\ldots,x-y_i,v)}+R,$$
(with $x-y_i$ occurring $k-1$ times), and 
$|R|\le\frac{2|x-y_i|^{m-1}|v|}{m!}\le\frac{2\delta^{m-1}}{m!}$.
As above, this implies
\begin{equation*}\begin{split}
|D(f-g)(x)|&\le\sum_{k=1}^{m-1}\frac{\delta^{k-1}}{k!}\|D^k(f-g)(y_i)\|_{op}+
\delta^{m-1}<5\delta^{m-1},
\end{split}\end{equation*}
and thus $\|f-g\|_1\le \epsilon$ if $\delta_1=\left(\frac{\epsilon}{5}
\right)^{\frac{1}{m-1}}.$

To bound the size of an $\epsilon$-net for $C^m_1(\X)$, it now only remains 
to count the number of possible arrays $A_f$ for $f\in C^m_1(\X)$.
Begin by counting the number of possibilities for the first column.  Since
$D^kf(y_1)$ is approximated in the $k$-1 entry of $A_f$ by a point from a 
$\frac{\delta^{m-k}}{2}$-net, the size of such a net is needed.  The 
space of symmetric $k$-linear forms is a finite-dimensional normed space,
and the size of a net for the unit ball of such a space is given in
Milman and Schechtman \cite{MS}, in terms of the dimension of the space.
To define an element $T$ of this space, it suffices to define $T(e_1,\ldots,
e_1,\ldots,e_d,\ldots,e_d),$ where $e_j$ appears $k_j$ times with $k_j\ge0$
for each $j$ and $\sum_{j=1}^dk_j=k$.  The number of such vectors $(k_j)$
is well-known (see, e.g., \cite{Ros}) to be $\binom{k+d-1}{k}$.  It follows
from the bound in \cite{MS} that there is a $\frac{\delta^{m-k}}{
2}$-net of the space of symmetric $k$-linear forms of size not greater than
$\left(1+\frac{4}{\delta^{m-k}}\right)^{\binom{k+d-1}{k}}\le \left(\frac{
5}{\delta^{m-k}}\right)^{\binom{k+d-1}{k}},$ assuming that $\delta<1$.
Since the only interesting case is $\epsilon\le 2$ (since $\|f-g\|_i\le 2$
for $i=0,1$ automatically), this is no restriction.
  The number of
possibilities for the first column of $A_f$ is thus bounded by
\begin{equation*}\begin{split}
\prod_{k=0}^{m-1}\left(\frac{5}{\delta^{m-k}}\right)^{\binom{k+d-1}{k}}
&\le\prod_{k=0}^{m-1}\left(\frac{5}{\delta^{m-k}}\right)^{\frac{(m+d-2)^k}{
k!}}\\&=\left(\frac{5}{\delta^m}\right)^{\sum_{k=0}^{m-1}\frac{(m+d
-2)^k}{k!}}\big(\delta\big)^{\sum_{k=0}^{m-1}k\frac{(m+d-2)^k}{k!}}\\&\le
\left(\frac{5}{\delta^m}\right)^{e^{m+d-2}}\big(\delta\big)^{m+d-2}\\&
=\exp\left[\left(\log(5)-m\log(\delta)\right)e^{m+d-2}+(m+d-2)\log(
\delta)\right]\\&\le\exp\left[\left(\log(5)-m\log(\delta)\right)e^{m+d-2}
\right],
\end{split}\end{equation*}
since $\delta<1$ and $m+d-2\ge 0$.
To bound the number of possibilities in the remaining columns, assume that
the $y_i$ have been ordered such that for all $j>1$, there is an $i<j$ 
with $|y_i-y_j|<2\delta$.  Now, for unit vectors $v_1,\ldots,v_k\subseteq\R^d$,
define the function $F(x):=\inprod{D^kf(x)}{(v_1,\ldots,v_k)},$ where the
dependence of $F$ on the $v_i$ has been suppressed.  By Taylor's theorem,
$$F(y_j)=\sum_{\ell=0}^{m-1-k}\frac{1}{\ell!}\inprod{D^{\ell}F(y_i)}{
(y_j-y_i,\ldots,y_j-y_i)}+R,$$
with $|R|\le\frac{(2\delta)^{m-k}}{(m-k)!}$.  
Let $A_f(i,j)$ denote the $i$-$j$-th entry
of the array $A_f$.  Then
\begin{equation*}\begin{split}
&\left|F(y_j)-\sum_{\ell=0}^{m-1-k}\frac{1}{\ell!}\inprod{A_f(i,k+\ell)}{(y_j-
y_i,\ldots,y_j-y_i,v_1,\ldots,v_k)}\right|\\&\qquad\le(2\delta)^{m-k}+
\left|\sum_{\ell=0}^{m-1-k}\frac{1}{\ell!}\inprod{D^{\ell}F(y_i)}{(y_j-
y_i,\ldots,y_j-y_i)}\right.\\&\qquad\qquad\qquad\qquad
\left.-\sum_{\ell=0}^{m-1-k}\frac{1}{\ell!}\inprod{A_f(i,k+
\ell)}{(y_j-y_i,\ldots,y_j-y_i,v_1,\ldots,v_k)}\right|\\&\qquad=
(2\delta)^{m-k}+\left|\sum_{\ell=0}^{m-1-k}\frac{1}{\ell!}\inprod{D^{\ell+k}
f(y_i)}{(y_j-y_i,\ldots,y_j-y_i,v_1,\ldots,v_k)}\right.\\&
\qquad\qquad\qquad\qquad
\left.-\sum_{\ell=0}^{m-1-k}\frac{1}{\ell!}\inprod{A_f(i,k+
\ell)}{(y_j-y_i,\ldots,y_j-y_i,v_1,\ldots,v_k)}\right|\\&\qquad\le
(2\delta)^{m-k}+\sum_{\ell=0}^{m-1-k}\frac{1}{\ell!}\left\|D^{k+\ell}f(y_i)
-A_f(i,k+\ell)\right\|_{op}(2\delta)^{\ell}\\&\qquad\le(2\delta)^{m-k}\left(
1+\frac{e}{2}\right).
\end{split}\end{equation*}
That is, given the information in the previous columns, the symmetric 
$k$-linear form $T(v_1,\ldots,v_k):=\inprod{D^kf(y_i)}{v_1,\ldots,v_k)}$ is
within a ball of radius $(2\delta)^{m-k}\left(1+\frac{e}{2}\right)$
with respect to the operator norm.  By the same argument that bounds the size
of the original $\frac{\delta^{m-k}}{2}$-net in the space, the number of
points of the net within this ball of radius $(2\delta)^{m-k}\left(1+
\frac{e}{2}\right)$ is bounded by 
$$\left(1+\frac{4(2\delta)^{m-k}\left(1+\frac{e}{2}\right)}{\delta^{m-
k}}\right)^{\binom{k+d-1}{k}}=\left(1+2^{m-k+2}\left(1+\frac{e}{2}\right)
\right)^{\binom{k+d-1}{k}}\le\left(2^{m-k+4}
\right)^{\binom{k+d-1}{k}} .$$
It follows that the number of possibilities for the column entries of $A_f$
after the first column is specified is bounded by 
\begin{equation*}\begin{split}
&\left[\prod_{k=0}^{m-1}\left(2^{m-k+4}
\right)^{\binom{k+d-1}{k}}\right]^{\frac{\vol(\X_1)}{\delta^d}}\\&\qquad=\exp\left[
\frac{\vol(\X_1)}{\delta^d}\sum_{k=0}^{m-1}\binom{k+d-1}{k}
\big((m-k+4)\log(2)\big)\right]\\&\qquad\le\exp\left[
\frac{\vol(\X_1)\big((m+4)\log(2)\big)}{\delta^d}\sum_{k=0}^{m-1}
\frac{1}{k!}(d+m-2)^k\right]\\&\qquad\le\exp\left[
\frac{\vol(\X_1)\big((m+4)\log(2)\big)}{\delta^d}e^{d+m-2}\right].
\end{split}\end{equation*}
It now follows that the total number of possible entries of $A_f$ is bounded
by 
$$\exp\left[\left(\log(5)-m
\log(\delta)+\frac{\vol(\X_1)\big((m
+4)\log(2)\big)}{\delta^d}\right)e^{d+m-2}\right].$$
Recall that $\delta_0$ and $\delta_1$ were chosen such that 
$\delta_0=\left(\frac{\epsilon}{5}
\right)^{\frac{1}{m}}$ and $\delta_1=\left(\frac{\epsilon}{5}
\right)^{\frac{1}{m-1}}$.
The $\epsilon$-covering number (for $\epsilon<2$)
of $C^m_1(\X)$ with respect to $\|\cdot\|_\infty$ is thus
bounded by
$$\exp\left[\left(2\log(5)-\log(\epsilon)+
\frac{c_0}{\epsilon^{\frac{d}{m}}}\right)e^{d+m-2}\right]$$
with
$$c_0=\vol(\X_1)\big((m+4)\log(2)\big)\big(5\big)^{\frac{d}{m}}.$$
The $\epsilon$-covering number
of $C^m_1(\X)$ for $\epsilon<2$ and $m\ge 2$ with respect to $\|\cdot\|_1$ is
bounded by 
$$\exp\left[\left((3\log(5)-\frac{m}{m-1}\log(\epsilon)+
\frac{c_1}{\epsilon^{\frac{d}{m-1}}}\right)e^{d+m-2}\right]$$
with
$$c_1=\vol(\X_1)\big((m+4)\log(2)\big)\big(5\big)^{\frac{d}{m-1}}.$$

\bigskip

\noindent {\bf Acknowledgements:} Many thanks to Michel Ledoux, Mark Meckes,
and Luke Rogers for helpful comments and suggestions.

\end{document}